\RequirePackage{ifpdf}
\ifpdf 
\documentclass[pdftex]{sigma}
\else
\documentclass{sigma}
\fi

\numberwithin{equation}{section}
\numberwithin{theorem}{section}
\numberwithin{proposition}{section}
\numberwithin{lemma}{section}
\numberwithin{corollary}{section}

\begin{document}

\allowdisplaybreaks

\renewcommand{\PaperNumber}{030}

\FirstPageHeading

\ShortArticleName{Rational Solutions of the Sasano System of Type $A_5^{(2)}$}

\ArticleName{Rational Solutions of the Sasano System of Type $\boldsymbol{A_5^{(2)}}$}

\Author{Kazuhide MATSUDA}

\AuthorNameForHeading{K.~Matsuda}

\Address{Department of Engineering Science, Niihama National College of Technology,\\
7-1 Yagumo-chou, Niihama, Ehime, 792-8580, Japan}
\Email{\href{mailto:matsuda@sci.niihama-nct.ac.jp}{matsuda@sci.niihama-nct.ac.jp}}

\ArticleDates{Received November 5, 2010, in f\/inal form March 17, 2011;  Published online March 25, 2011}

\Abstract{In this paper,
we completely classify the rational solutions of the Sasano system of type $A_5^{(2)}$,
which is given by the coupled Painlev\'e III system.
This system of dif\/ferential equations has the af\/f\/ine Weyl group symmetry of type $A_5^{(2)}$.}


\Classification{33E17; 34M55}

\section{Introduction}

Paul Painlev\'e and his colleagues
\cite{Painleve, Gambier}
intended to f\/ind new transcendental functions def\/ined by second order nonlinear dif\/ferential equations.
In general,
nonlinear dif\/ferential equations have moving branch points.
If a solution has moving branch points, it is too complicated and is not worth considering.
Therefore,
they determined the second order nonlinear dif\/ferential equations with rational coef\/f\/icients
which have no moving branch points.
As a result,
the standard forms of such equations turned out to be given
by the following six equations:
\begin{alignat*}{3}
&P_{\rm I} :   \qquad  &&
 y^{\prime \prime}
=
6
y^2+t, &  \\
&P_{\rm II} : \qquad  &&
y^{\prime \prime}
=
2 y^3 +ty+\alpha, &  \\
&P_{\rm III} :  \quad &&
y^{\prime \prime}
=
\frac{1}{y}(y^{\prime})^2
-\frac{1}{t}y^{\prime}
+\frac{1}{t}(\alpha y^2+\beta)
+\gamma y^3 + \frac{\delta}{y}, &  \\
&P_{\rm IV} : \qquad &&
y^{\prime \prime}
=
\frac{1}{2y}(y^{\prime})^2
+\frac32 y^3 +4ty^2
+2(t^2-\alpha)y+\frac{\beta}{y}, &  \\
&P_{\rm V} :  \qquad  &&
y^{\prime \prime}
=
\left(
\frac{1}{2y}+\frac{1}{y-1}
\right)
(y^{\prime})^2
-\frac{1}{t}y^{\prime}
+\frac{(y-1)^2}{t^2} \left(\alpha y + \frac{\beta}{y} \right)
+\frac{\gamma}{t}y
+\delta \frac{y(y+1)}{y-1},&  \\
&P_{\rm VI}  : \quad   &&
y^{\prime \prime}
=
\frac12
\left(
\frac{1}{y}+\frac{1}{y-1}+\frac{1}{y-t}
\right)
(y^{\prime})^2
-\left(\frac{1}{t}+\frac{1}{t-1}+\frac{1}{y-t}\right)y^{\prime} & \\
&  & &
{} \phantom{y^{\prime \prime}=}{} +\frac{y(y-1)(y-t)}{t^2(t-1)^2}
\left(
\alpha+\beta \frac{t}{y^2}+\gamma \frac{t-1}{(y-1)^2}
+\delta \frac{t(t-1)}{(y-t)^2}
\right),&
\end{alignat*}
where $^\prime=d/dt$
and
$\alpha$, $\beta$, $\gamma$, $\delta$ are all complex parameters.
In this article,
our concern is with the B\"acklund transformations and special solutions which are given by rational, algebraic functions or classical special functions.

Each of $P_J$ $(J={\rm II}, {\rm III}, {\rm IV}, {\rm V}, {\rm VI})$ has B\"acklund transformations,
which transform solutions into other solutions of the same equation
with dif\/ferent parameters.
It was shown by Okamoto \cite{oka1,oka2,oka3,oka4} that
the B\"acklund transformation groups of the Painlev\'e equations except for $P_{\rm I}$ are isomorphic
to the extended af\/f\/ine Weyl groups.
For $ P_{\rm II}$, $P_{\rm III}$, $P_{\rm IV}$, $P_{\rm V}$, and $P_{\rm VI}$,
the B\"acklund transformation groups correspond to
$A^{ ( 1 ) }_1$,
$A^{ ( 1 ) }_1 \bigoplus A^{ ( 1 ) }_1$,
$A^{ ( 1 ) }_2$,
$A^{ ( 3 ) }_3$,
and
$D^{ ( 1 ) }_4$,
respectively.

While generic solutions of the Painlev\'e equations are
``new transcendental functions'',
there are special solutions which are expressible
in terms of rational, algebraic, or classical special functions.

For example,
Airault \cite{Airault} constructed explicit rational solutions of $P_{\rm II}$ and $P_{\rm IV}$
with their B\"acklund transformations.
Milne, Clarkson and Bassom \cite{Milne} treated $P_{\rm III}$,  and
described their B\"acklund transformations
and exact solution hierarchies, which are given by rational, algebraic, or certain Bessel functions.
Bassom, Clarkson and Hicks \cite{Bassom-Clarkson-Hicks}
dealt with $P_{\rm IV}$, and
described their B\"acklund transformations and
exact solution hierarchies, which are expressed by rational functions, the parabolic cylinder functions or the complementary error functions.
Clarkson~\cite{Clarkson}
studied some rational and algebraic solutions of $P_{\rm III}$ and
showed that
these solutions are expressible
in terms of special polynomials def\/ined by second order, bilinear dif\/ferential-dif\/ference equations
which are equivalent to Toda equations.

Furthermore,
the rational solutions of $P_J$ $(J={\rm II},{\rm III},{\rm IV},{\rm V},{\rm VI})$
were classif\/ied by
Yablonski and Vorobev~\cite{Yab:59,Vorob},
Gromak~\cite{Gr:83,Gro},
Murata~\cite{Mura1, Mura2},
Kitaev, Law and McLeod~\cite{Kit-Law-McL},
Mazzoco~\cite{Mazzo1} and Yuang and Li~\cite{YuangLi}.

Noumi and Yamada \cite{NoumiYamada-B} discovered the equation of type $A^{(1)}_l$ $(l\geq 2)$,
whose B\"acklund
transformation group is isomorphic to the extended af\/f\/ine Weyl group $\tilde{W}(A^{(1)}_l)$.
The Noumi and Yamada systems of types $A_2^{(1)}$ and $A_3^{(1)}$ correspond to the fourth and f\/ifth
Painlev\'e equations, respectively.
Moreover,
we \cite{Matsuda1, Matsuda2} classif\/ied the rational solutions of the Noumi and Yamada systems of types~$A_4^{(1)}$ and~$A_5^{(1)}$.

Sasano \cite{Sasano-1} found
the coupled Painlev\'e V and VI systems
which have
the af\/f\/ine Weyl group symmetries of types $D^{(1)}_5$ and $D_6^{(1)}$.
In addition, he \cite{Sasano-2}
obtained the equation of the af\/f\/ine
Weyl group symmetry of type $A^{(2)}_5$,
which is def\/ined by
\begin{gather*}
A_5^{(2)}(\alpha_j)_{0\leq j \leq 3} \
\begin{cases}
tq_1^{\prime} =2q_1^2p_1-q_1^2+(\alpha_0+\alpha_1+\alpha_3)q_1-t+4tp_2+2q_1q_2p_2,  \\
tp_1^{\prime} =-2q_1p_1^2+2q_1p_1-(\alpha_0+\alpha_1+\alpha_3)p_1+\alpha_0-2p_1q_2p_2,  \\
tq_2^{\prime} =2q_2^2p_2-q_2^2+(\alpha_0+\alpha_1+\alpha_3)q_2-t+4tp_1+2q_1p_1q_2, \\
tp_2^{\prime} =-2q_2p_2^2+2q_2p_2-(\alpha_0+\alpha_1+\alpha_3)p_2+\alpha_1-2q_1p_1p_2, \\
\alpha_0+\alpha_1+2\alpha_2+\alpha_3=1/2,
\end{cases}
\end{gather*}
where $^\prime=d/dt$.
This system of dif\/ferential equations is also expressed by the Hamiltonian system:
\begin{gather*}
t\frac{dq_1}{dt}=\frac{\partial H}{\partial p_1},   \qquad
t\frac{dp_1}{dt}=-\frac{\partial H}{\partial q_1},   \qquad
t\frac{dq_2}{dt}=\frac{\partial H}{\partial p_2},   \qquad
t\frac{dp_2}{dt}=-\frac{\partial H}{\partial q_2},
\end{gather*}
where
the Hamiltonian $H$ is given by
\begin{gather*}
H=
q_1^2p_1^2-q_1^2p_1+(\alpha_0+\alpha_1+\alpha_3)q_1p_1-\alpha_0q_1-tp_1  \\
\phantom{H=}{}
+q_2^2p_2^2-q_2^2p_2+(\alpha_0+\alpha_1+\alpha_3)q_2p_2-\alpha_1q_2-tp_2+4tp_1p_2+2q_1p_1q_2p_2.
\end{gather*}
Let us note that
Mazzocco and Mo \cite{Mazzo2} studied the Hamiltonian structure of the $P_{\rm II}$ hierarchy,
and
Hone \cite{Hone} studied the coupled Painlev\'e systems from the similarity reduction of the
Hirota--Satsuma system and another gauge-related system,
and
presented their B\"acklund transformations and special solutions.

$A_5^{(2)}(\alpha_j)_{0\leq j \leq 3}$
has the B\"acklund transformations $s_0$, $s_1$, $s_2$, $s_3$, $\pi$,
which are given by
\begin{alignat*}{5}
&s_0:\quad & &(*) & &\rightarrow & &\left(
q_1+\frac{\alpha_0}{p_1},p_1,q_2,p_2,t;-\alpha_0,\alpha_1,\alpha_2+\alpha_0,\alpha_3\right), & \\
&s_1: & & (*) & &\rightarrow &&\left(
q_1,p_1,q_2+\frac{\alpha_1}{p_2},p_2,t;\alpha_0,-\alpha_1,\alpha_2+\alpha_1,\alpha_3
\right), &  \\
&s_2: & &(*) &
&\rightarrow &
&\left(
q_1, p_1-\frac{\alpha_2q_2}{q_1q_2+t}, q_2, p_2+\frac{\alpha_2q_1}{q_1q_2+t},t; \alpha_0+\alpha_2,\alpha_1+\alpha_2,-\alpha_2,\alpha_3+2\alpha_2
\right), &  \\
&s_3: & &(*) &
&\rightarrow &
&\left(
q_1+\frac{\alpha_3}{p_1+p_2-1},p_1,q_2+\frac{\alpha_3}{p_1+p_2-1}, p_2,t;\alpha_0,\alpha_1,\alpha_2+\alpha_3,-\alpha_3
\right),&  \\
&\pi:  & &(*) &
&\rightarrow &
&\left(
q_2,p_2,q_1,p_1,t;\alpha_1,\alpha_0,\alpha_2,\alpha_3
\right), &
\end{alignat*}
with the notation $(*)=(q_1,p_1,q_2,p_2,t;\alpha_0,\alpha_1,\alpha_2,\alpha_3)$.
The B\"acklund transformation group $\langle s_0, s_1, s_2, s_3, \pi \rangle$
is isomorphic to the af\/f\/ine Weyl group of type $A_5^{(2)}$.

Our main theorem is as follows:

\begin{theorem}
For a rational solution of $A_5^{(2)}(\alpha_j)_{0\leq j \leq 3}$,
by some B\"acklund transformations,
the solution and parameters
can be transformed so that
\begin{gather*}
(q_1,p_1,q_2,p_2)=(0,1/4,0,1/4) \qquad {\it and} \\
 (\alpha_0,\alpha_1,\alpha_2,\alpha_3)=(\alpha_3/2,\alpha_3/2,\alpha_2,\alpha_3)=(\alpha_3/2,\alpha_3/2, 1/4-\alpha_3,\alpha_3),
\end{gather*}
respectively.
Furthermore,
for $A_5^{(2)}(\alpha_j)_{0\leq j \leq 3}$,
there exists a rational solution
if and only if
one of the following occurs:
\begin{alignat*}{3}
&(1) &\quad  -2\alpha_0+\alpha_3&\in\mathbb{Z}, &\quad -2\alpha_1+\alpha_3&\in\mathbb{Z},  \\
&(2) &  -2\alpha_0+\alpha_3&\in\mathbb{Z}, &  2\alpha_1+\alpha_3&\in\mathbb{Z},   \\
&(3) &   2\alpha_0+\alpha_3&\in\mathbb{Z},  & -2\alpha_1+\alpha_3&\in\mathbb{Z},  \\
&(4) &   2\alpha_0+\alpha_3&\in\mathbb{Z},  &  2\alpha_1+\alpha_3&\in\mathbb{Z},  \\
&(5) &  -2\alpha_0+\alpha_3&\in\mathbb{Z},  & \alpha_3-1/2&\in\mathbb{Z},\\
&(6) &  -2\alpha_1+\alpha_3&\in\mathbb{Z},  & \alpha_3-1/2&\in\mathbb{Z}.
\end{alignat*}
\end{theorem}

This paper is organized as follows.
In Section~\ref{section1}, for $A_5^{(2)}(\alpha_j)_{0\leq j \leq 3}$,
we determine meromorphic solutions at $t=\infty$.
Then,
we f\/ind that
the constant terms $a_{\infty,0}$, $c_{\infty,0}$ of the Laurent series of $q_1$, $q_2$ at $t=\infty$ are
given by
\[
a_{\infty,0}:=-2\alpha_0+\alpha_3, \qquad c_{\infty,0}:=-2\alpha_1+\alpha_3,
\]
respectively.

In Section~\ref{section2},
for $A_5^{(2)}(\alpha_j)_{0\leq j \leq 3}$,
we determine meromorphic solutions at $t=0$.
Then,
we see that
the constant terms $a_{0,0}$, $c_{0,0}$ of the Laurent series of $q_1$, $q_2$ at $t=0$ are
given by the parameters $\alpha_0$, $\alpha_1$, $\alpha_2$, $\alpha_3$.

In Section~\ref{section3},
for $A_5^{(2)}(\alpha_j)_{0\leq j \leq 3}$,
we treat meromorphic solutions at $t=c\in\mathbb{C}^{\ast}$,
where  in this paper,
$\mathbb{C}^{*}$ means the set of nonzero complex numbers.
Then,
we observe that
$q_1$, $q_2$ have both a pole of order of at most one at $t=c$
and the residues of $q_1$, $q_2$ at $t=c$ are expressed by $nc$ $(n\in\mathbb{Z})$.
Thus,
it follows that
\begin{equation}
\label{eqn:formula-giving-necessary condition}
a_{\infty,0}-a_{0,0} \in\mathbb{Z}, \qquad c_{\infty,0}-c_{0,0} \in\mathbb{Z},
\end{equation}
which gives a necessary condition for $A_5^{(2)}(\alpha_j)_{0\leq j \leq 3}$
to have a rational solution.

In Section~\ref{section4},
using the meromorphic solution at $t=\infty, 0$,
we f\/irst compute the constant terms of the Laurent series of the Hamiltonian at $t=\infty,0$.
Furthermore,
by the meromorphic solution at $=c\in\mathbb{C}^{*}$,
we calculate the residue of $H$ at $t=c$.

In Section~\ref{section5},
by equation (\ref{eqn:formula-giving-necessary condition}),
we obtain the necessary conditions for $A_5^{(2)}(\alpha_j)_{0\leq j \leq 3}$
to have rational solutions,
which are given in our main theorem.
Furthermore,
we show that
if there exists a rational solution for $A_5^{(2)}(\alpha_j)_{0\leq j \leq 3}$,
the parameters can be transformed so that $-2\alpha_0+\alpha_3\in\mathbb{Z}$, $-2\alpha_1+\alpha_3\in\mathbb{Z}$.

In Section~\ref{section6},
we def\/ine shift operators,
and
for a rational solution of $A_5^{(2)}(\alpha_j)_{0\leq j \leq 3}$,
we transform the parameters to
\[
(\alpha_0,\alpha_1,\alpha_2,\alpha_3)=(\alpha_3/2,\alpha_3/2,\alpha_2,\alpha_3).
\]

In Section~\ref{section7},
we determine rational solutions of $A_5^{(2)}(\alpha_3/2,\alpha_3/2,\alpha_2,\alpha_3)$
and
prove our main theorem.

In Appendix~\ref{appendixA},
using the shift operators,
we give examples of rational solutions.

\section[Meromorphic solutions at $t=\infty$]{Meromorphic solutions at $\boldsymbol{t=\infty}$}\label{section1}

In this section, for $A_5^{(2)}(\alpha_j)_{0\leq j \leq 3}$,
we treat meromorphic solutions at $t=\infty$.
For the purpose,
in this paper,
we def\/ine the coef\/f\/icients of the Laurent series of $q_1$, $p_1$, $q_2$, $p_2$ at $t=\infty$
by
$a_{\infty,k}$, $b_{\infty,k}$, $c_{\infty,k}$, $d_{\infty,k}$, $k\in\mathbb{Z}$.

\subsection[The case where $q_1$, $p_1$, $q_2$, $p_2$ are all holomorphic at $t=\infty$]{The case where $\boldsymbol{q_1}$, $\boldsymbol{p_1}$, $\boldsymbol{q_2}$, $\boldsymbol{p_2}$ are all holomorphic at $\boldsymbol{t=\infty}$}

\begin{proposition}
Suppose that
for $A_5^{(2)}(\alpha_j)_{0\leq j \leq 3}$,
there exists a solution such that
$q_1$, $p_1$, $q_2$, $p_2$ are all holomorphic at $t=\infty$.
Then,
\begin{gather*}
\begin{cases}
q_1=(-2\alpha_0+\alpha_3)+\cdots, \\
p_1=1/4+(-2\alpha_1+\alpha_3)(-2\alpha_1-\alpha_3) t^{-1}/4+\cdots, \\
q_2=(-2\alpha_1+\alpha_3)+\cdots, \\
p_2=1/4+(-2\alpha_0+\alpha_3)(-2\alpha_0-\alpha_3)t^{-1}/4 +\cdots.
\end{cases}
\end{gather*}
\end{proposition}


\begin{proposition}
\label{prop:uniqueness}
Suppose that
for $A_5^{(2)}(\alpha_j)_{0\leq j \leq 3}$,
there exists a solution such that
$q_1$, $p_1$, $q_2$, $p_2$ are all holomorphic at $t=\infty$.
Then,
it is unique.
\end{proposition}

\begin{proof}
We set
\begin{gather*}
q_1  =a_{\infty,0} +a_{\infty,-1}t^{-1}+\cdots +a_{\infty,-(k-1)}t^{-(k-1)}+a_{\infty,-k}t^{-k}+a_{\infty,-(k+1)}t^{-(k+1)}+\cdots,   \\
p_1  =1/4 +b_{\infty,-1}t^{-1}+\cdots +b_{\infty,-(k-1)}t^{-(k-1)}+b_{\infty,-k}t^{-k}+b_{\infty,-(k+1)}t^{-(k+1)}+\cdots,   \\
q_2  =c_{\infty,0} +c_{\infty,-1}t^{-1}+\cdots +c_{\infty,-(k-1)}t^{-(k-1)}+c_{\infty,-k}t^{-k}+c_{\infty,-(k+1)}t^{-(k+1)}+\cdots,   \\
p_2  =1/4  +d_{\infty,-1}t^{-1}+\cdots +d_{\infty,-(k-1)}t^{-(k-1)}+d_{\infty,-k}t^{-k}+d_{\infty,-(k+1)}t^{-(k+1)}+\cdots,
\end{gather*}
where $a_{\infty,0}$,  $b_{\infty,-1}$, $c_{\infty,0}$, $d_{\infty,-1}$ all have been determined.
\par
Comparing the coef\/f\/icients of the terms $t^{-k}$ $(k \geq 1)$ in
\begin{gather*}
tp_1^{\prime} =-2q_1p_1^2+2q_1p_1-(\alpha_0+\alpha_1+\alpha_3)p_1+\alpha_0-2p_1q_2p_2,  \\
tp_2^{\prime} =-2q_2p_2^2+2q_2p_2-(\alpha_0+\alpha_1+\alpha_3)p_2+\alpha_1-2q_1p_1p_2,
\end{gather*}
we have
\begin{gather*}
3a_{\infty,-k}/8-c_{\infty,-k}/8
 =-kb_{\infty,-k}+(\alpha_0+\alpha_1+\alpha_3)b_{\infty,-k}  \\
  \qquad {}+2\sum a_{\infty,-l} b_{\infty,-m} b_{\infty,-n}
-2\sum a_{\infty,-l} b_{\infty,-m}
+2\sum c_{\infty,-l} b_{\infty,-m} d_{\infty,-n},   \\
-a_{\infty,-k}/8+3c_{\infty,-k}/8
 =-kd_{\infty,-k} +(\alpha_0+\alpha_1+\alpha_3)d_{\infty,-k}  \\
  \qquad {}+2\sum c_{\infty,-l} d_{\infty,-m} d_{\infty,-n}
-2\sum c_{\infty,-l} d_{\infty,-m}
+2\sum a_{\infty,-l} b_{\infty,-m} d_{\infty,-n},
\end{gather*}
where the f\/irst and third sums extend over nonnegative integers $l$, $m$, $n$ such that $l+m+n=k$ and $0\leq l <k$,
and the second sums extend over nonnegative integers $l$, $m$ such that $l+m=k$ and $m\geq 1$.
Therefore,
$a_{\infty,-k}$, $c_{\infty,-k}$ are both inductively determined.

Comparing the coef\/f\/icients of the terms $t^{-k}$ $(k \geq 1)$ in
\begin{gather*}
tq_1^{\prime} =2q_1^2p_1-q_1^2+(\alpha_0+\alpha_1+\alpha_3)q_1-t+4tp_2+2q_1q_2p_2,  \\
tq_2^{\prime} =2q_2^2p_2-q_2^2+(\alpha_0+\alpha_1+\alpha_3)q_2-t+4tp_1+2q_1p_1q_2,
\end{gather*}
we obtain
\begin{gather*}
4d_{\infty,-(k+1)}
 =-ka_{\infty,-k}-(\alpha_0+\alpha_1+\alpha_3)a_{\infty,-k}  \\
 \qquad {}-2\sum a_{\infty,-l}a_{\infty,-m}b_{\infty,-n}+\sum a_{\infty,-l}a_{\infty,-m}
-2\sum a_{\infty,-l}c_{\infty,-m}d_{\infty,-n}, \\
4b_{\infty,-(k+1)}
 =-kc_{\infty,-k}-(\alpha_0+\alpha_1+\alpha_3)c_{\infty,-k}  \\
 \qquad {} -2\sum c_{\infty,-l}c_{\infty,-m}d_{\infty,-n}+\sum c_{\infty,-l}c_{\infty,-m}
-2\sum c_{\infty,-l}a_{\infty,-m}b_{\infty,-n},
\end{gather*}
where the f\/irst and third sums extend over nonnegative integers $l$, $m$, $n$ such that $l+m+n=k$,
and the second sums extend over nonnegative integers $l$, $m$ such that $l+m=k$.
Therefore,
$b_{\infty,-(k+1)}$, $d_{\infty,-(k+1)}$ are both inductively determined,
which proves the proposition.
\end{proof}

\subsection[The case where one of $(q_1,p_1,q_2,p_2)$ has a pole at $t=\infty$]{The case where one of $\boldsymbol{(q_1,p_1,q_2,p_2)}$ has a pole at $\boldsymbol{t=\infty}$}

In this subsection,
we deal with the case in which
one of $(q_1,p_1,q_2,p_2)$ has a pole at $t=\infty$.
For the purpose,
by $\pi$,
we have only to consider the following two cases:
\begin{enumerate}\itemsep=0pt
\item[(1)] $q_1$ has a pole at $t=\infty$ and $p_1$, $q_2$, $p_2$ are all holomorphic at $t=\infty$,

\item[(2)] $p_1$ has a pole at $t=\infty$ and $q_1$, $q_2$, $p_2$  are all holomorphic at $t=\infty$.
\end{enumerate}

\subsubsection[The case where $q_1$ has a pole at $t=\infty$]{The case where $\boldsymbol{q_1}$ has a pole at $\boldsymbol{t=\infty}$}

\begin{proposition}
For $A_5^{(2)}(\alpha_j)_{0\leq j \leq 3}$,
there exists no solution such that
$q_1$ has a pole at $t=\infty$ and $p_1$, $q_2$, $p_2$ are all holomorphic at $t=\infty$.
\end{proposition}


\subsubsection[The case where $p_1$ has a pole at $t=\infty$]{The case where $\boldsymbol{p_1}$ has a pole at $\boldsymbol{t=\infty}$}

\begin{proposition}
For $A_5^{(2)}(\alpha_j)_{0\leq j \leq 3}$,
there exists no solution such that
$p_1$ has a pole at $t=\infty$ and $q_1$, $q_2$, $p_2$ are all holomorphic at $t=\infty$.
\end{proposition}


\subsection[The case where two of $(q_1,p_1,q_2,p_2)$ have a pole at $t=\infty$]{The case where two of $\boldsymbol{(q_1,p_1,q_2,p_2)}$ have a pole at $\boldsymbol{t=\infty}$}

In this subsection,
we deal with the case in which
two of $(q_1,p_1,q_2,p_2)$ has a pole at $t=\infty$.
For the purpose,
by $\pi$,
we have only to consider the following four cases:
\begin{enumerate}\itemsep=0pt
\item[(1)]  $q_1$, $p_1$ have both a pole at $t=\infty$ and $q_2$, $p_2$ are both holomorphic at $t=\infty$,
\item[(2)] $q_1$, $q_2$ have both a pole at $t=\infty$ and $p_1$, $p_2$ are both holomorphic at $t=\infty$,
\item[(3)] $q_1$, $p_2$ have both a pole at $t=\infty$ and $p_1$, $q_2$ are both holomorphic at $t=\infty$,
\item[(4)] $p_1$, $p_2$ have both a pole at $t=\infty$ and $q_1$, $q_2$ are both holomorphic at $t=\infty$.
\end{enumerate}

\subsubsection[The case where $q_1$, $p_1$ have a pole at $t=\infty$]{The case where $\boldsymbol{q_1}$, $\boldsymbol{p_1}$ have a pole at $\boldsymbol{t=\infty}$}

\begin{proposition}
For $A_5^{(2)}(\alpha_j)_{0\leq j \leq 3}$,
there exists no solution such that
$q_1$, $p_1$ have both a pole at $t=\infty$ and $q_2$, $p_2$ are both holomorphic at $t=\infty$.
\end{proposition}


\subsubsection[The case where $q_1$, $q_2$ have a pole at $t=\infty$]{The case where $\boldsymbol{q_1}$, $\boldsymbol{q_2}$ have a pole at $\boldsymbol{t=\infty}$}

\begin{proposition}
For $A_5^{(2)}(\alpha_j)_{0\leq j \leq 3}$,
there exists no solution such that
$q_1$, $q_2$ have both a pole at $t=\infty$ and $p_1$, $p_2$ are both holomorphic at $t=\infty$.
\end{proposition}


\subsubsection[The case where $q_1$, $p_2$ have a pole at $t=\infty$]{The case where $\boldsymbol{q_1}$, $\boldsymbol{p_2}$ have a pole at $\boldsymbol{t=\infty}$}

By direct calculation,
we can obtain the following two lemmas:

\begin{lemma}
\label{lem:q_1=0}
Suppose that
for $A_5^{(2)}(\alpha_j)_{0\leq j \leq 3}$,
$q_1\equiv 0$.
Then, one of the following occurs:
\begin{alignat*}{3}
&(1)\quad &&
\alpha_0=\frac14, \qquad \alpha_3=\frac12,   \qquad       \mbox{and}      & \\
&&& (q_1,p_1,q_2,p_2)= \left(
0,   \frac14+\frac{(4\alpha_1-1)(4\alpha_1+1)}{16t},  -2\alpha_1+\frac12,  \frac14
\right), &   \\
&(2)\quad &&
\alpha_0=\frac{\alpha_3}{2}, \qquad \alpha_1=\frac{\alpha_3}{2},    \qquad \mbox{and}     \qquad (q_1,p_1,q_2,p_2)=
\left(
0,   \frac14,                                           0,        \frac14
\right), &      \\
&(3)\quad && \alpha_0=\frac{\alpha_3}{2}, \qquad \alpha_1=-\frac{\alpha_3}{2},     \qquad \mbox{and}  \qquad (q_1,p_1,q_2,p_2)=
\left(
0,   \frac14,                                        2\alpha_3,            \frac14
\right). &
\end{alignat*}
\end{lemma}

\begin{lemma}
\label{lem:q_2=0}
Suppose that
for $A_5^{(2)}(\alpha_j)_{0\leq j \leq 3}$,
$q_2\equiv 0$.
Then, one of the following occurs:
\begin{alignat*}{3}
&(1)\quad && \alpha_1=\frac14, \qquad \alpha_3=\frac12,      \qquad \mbox{and} &  \\
&&& (q_1,p_1,q_2,p_2) =
\left (
-2\alpha_0+\frac12,  \frac14,   0,   \frac14+\frac{(4\alpha_0-1)(4\alpha_0+1)}{16t}
\right),  &    \\
&(2)\quad && \alpha_0=\frac{\alpha_3}{2}, \qquad \alpha_1=\frac{\alpha_3}{2},     \qquad \mbox{and}     \qquad (q_1,p_1,q_2,p_2)=
\left(
 0,                        \frac14,   0,    \frac14
\right),  &      \\
&(3)\quad && \alpha_0=-\frac{\alpha_3}{2},   \qquad \alpha_1=\frac{\alpha_3}{2},   \qquad \mbox{and}     \qquad (q_1,p_1,q_2,p_2)=
\left(
 2\alpha_3,                               \frac14,    0,    \frac14
\right).&
\end{alignat*}
\end{lemma}

By Lemma~\ref{lem:q_2=0},
we f\/ind that $q_2\not\equiv 0$.
Now,
let us assume that
$q_1$ has a pole of order $n_0$ $(n_0\geq 1)$
and
$p_2$ has a pole of order $n_3$ $(n_3\geq 1)$.

\begin{lemma}
For $A_5^{(2)}(\alpha_j)_{0\leq j \leq 3}$,
there exists a solution such that
$q_1$, $p_2$ have both a pole at $t=\infty$ and $p_1$, $q_2$ are both holomorphic at $t=\infty$.
Then,
$n_0\neq n_3$.
\end{lemma}

\begin{proof}
We suppose that $n_0=n_3$.
Especially,
we treat the case where $n_0=n_3=1$ and show contradiction.
If $n_0=n_3>1$,
we can prove contradiction in the same way.

Comparing the coef\/f\/icients of the terms $t^2$, $t$
in
\begin{gather*}
tq_1^{\prime} =2q_1^2p_1-q_1^2+(\alpha_0+\alpha_1+\alpha_3)q_1-t+4tp_2+2q_1q_2p_2,  \\
tp_1^{\prime} =-2q_1p_1^2+2q_1p_1-(\alpha_0+\alpha_1+\alpha_3)p_1+\alpha_0-2p_1q_2p_2,
\end{gather*}
we have
\begin{gather}
2a_{\infty,1}^2b_{\infty,0}-a_{\infty,1}^2+4d_{\infty,1}+2a_{\infty,1}c_{\infty,0}d_{\infty,1}=0, \nonumber\\
-2a_{\infty,1}b_{\infty,0}^2+2a_{\infty,1}b_{\infty,0}-2b_{\infty,0}c_{\infty,0}d_{\infty,1}=0,
\label{*}
\end{gather}
respectively.

Comparing the coef\/f\/icients of the terms $t$, $t^2$
in
\begin{gather*}
tq_2^{\prime} =2q_2^2p_2-q_2^2+(\alpha_0+\alpha_1+\alpha_3)q_2-t+4tp_1+2q_1p_1q_2, \\
tp_2^{\prime} =-2q_2p_2^2+2q_2p_2-(\alpha_0+\alpha_1+\alpha_3)p_2+\alpha_1-2q_1p_1p_2,
\end{gather*}
we obtain
\begin{gather}
2c_{\infty,0}^2d_{\infty,1}-1+4b_{\infty,0}+2a_{\infty,1}b_{\infty,0}c_{\infty,0}=0,  \nonumber\\
-2c_{\infty,0}d_{\infty,1}^2-2a_{\infty,1}b_{\infty,0}d_{\infty,1}=0,
\label{**}
\end{gather}
which implies that $b_{\infty,0}=1/4$.
Furthermore,
from the second equation in \eqref{*} and the f\/irst equation in \eqref{**},
it follows that $a_{\infty,1}=0$,
which is impossible.
\end{proof}

\begin{lemma}
Suppose that
for $A_5^{(2)}(\alpha_j)_{0\leq j \leq 3}$,
there exists a solution such that
$q_1$, $p_2$ have both a pole at $t=\infty$
and
$p_1$, $q_2$ are both holomorphic at $t=\infty$.
Then, $n_0<n_3$.
\end{lemma}


\begin{proposition}
For $A_5^{(2)}(\alpha_j)_{0\leq j \leq 3}$,
there exists no solution such that
$q_1$, $p_2$ have both a pole at $t=\infty$ and $p_1$, $q_2$ are both holomorphic at $t=\infty$.
\end{proposition}

\begin{proof}
We treat the case where $(n_0,n_3)=(1,2)$ and show contradiction.
The other cases can be proved in the same way.

Comparing the coef\/f\/icients of the terms $t^3$
in
\begin{equation*}
tq_1^{\prime} =2q_1^2p_1-q_1^2+(\alpha_0+\alpha_1+\alpha_3)q_1-t+4tp_2+2q_1q_2p_2,
\end{equation*}
we have $d_{\infty,2}=0$,
which is impossible.
\end{proof}

\subsubsection[The case where $p_1$, $p_2$ have a pole at $t=\infty$]{The case where $\boldsymbol{p_1}$, $\boldsymbol{p_2}$ have a pole at $\boldsymbol{t=\infty}$}

\begin{proposition}
For $A_5^{(2)}(\alpha_j)_{0\leq j \leq 3}$,
there exists no solution such that
$p_1$, $p_2$ have both a pole at $t=\infty$ and $q_1$, $q_2$ are both holomorphic at $t=\infty$.
\end{proposition}


\subsection[The case where three of $(q_1,p_1,q_2,p_2)$ have a pole at $t=\infty$]{The case where three of $\boldsymbol{(q_1,p_1,q_2,p_2)}$ have a pole at $\boldsymbol{t=\infty}$}
In this subsection,
considering $\pi$,
we treat the following two cases:
\begin{enumerate}\itemsep=0pt
\item[(1)]  $q_1$, $p_1$, $q_2$ all have a pole at $t=\infty$ and $p_2$ is holomorphic at $t=\infty$,
\item[(2)] $q_1$, $p_1$, $p_2$ all have a pole at $t=\infty$ and $q_2$ is holomorphic at $t=\infty$.
\end{enumerate}

\subsubsection[The case where $q_1$, $p_1$, $q_2$ have a pole at $t=\infty$]{The case where $\boldsymbol{q_1}$, $\boldsymbol{p_1}$, $\boldsymbol{q_2}$ have a pole at $\boldsymbol{t=\infty}$}

\begin{proposition}
For $A_5^{(2)}(\alpha_j)_{0\leq j \leq 3}$,
there exists no solution such that
$q_1$, $p_1$, $q_2$ all have a~pole at $t=\infty$ and $p_2$ is holomorphic at $t=\infty$.
\end{proposition}


\subsubsection[The case where $q_1$, $p_1$, $p_2$ have a pole at $t=\infty$]{The case where $\boldsymbol{q_1}$, $\boldsymbol{p_1}$, $\boldsymbol{p_2}$ have a pole at $\boldsymbol{t=\infty}$}

By Lemma \ref{lem:q_2=0},
let us note that $q_2\not\equiv 0$.

\begin{lemma}
Suppose that
for $A_5^{(2)}(\alpha_j)_{0\leq j \leq 3}$,
there exists a solution such that
$q_1$, $p_1$, $p_2$ all have a pole at $t=\infty$ and $q_2$ is holomorphic at $t=\infty$.
Moreover,
assume that
$q_1$, $p_1$, $p_2$ has a pole of order $n_0$, $n_1$, $n_3$ $(n_0, n_1, n_3\geq 1)$ at $t=\infty$,
respectively.
Then,
$n_3\geq n_0+n_1$.
\end{lemma}

\begin{proof}
Considering that
\[
tp_1^{\prime} =-2q_1p_1^2+2q_1p_1-(\alpha_0+\alpha_1+\alpha_3)p_1+\alpha_0-2p_1q_2p_2,
\]
we can prove the lemma.
\end{proof}

Therefore,
we def\/ine the nonnegative integer $k$ by $n_3=n_0+n_1+k$.

\begin{lemma}
Suppose that
for $A_5^{(2)}(\alpha_j)_{0\leq j \leq 3}$,
there exists a solution such that
$q_1$, $p_1$, $p_2$ all have a pole at $t=\infty$ and $q_2$ is holomorphic at $t=\infty$.
Then,
$c_{\infty,0}=c_{\infty,-1}=\cdots=c_{\infty,-(k-1)}=0$, $a_{\infty,n_0}b_{\infty,1}+c_{\infty,-k}d_{\infty,n_3}=0$,
and
$n_0-k\geq 1$.
\end{lemma}

\begin{proof}
Considering that
\[
tp_1^{\prime} =-2q_1p_1^2+2q_1p_1-(\alpha_0+\alpha_1+\alpha_3)p_1+\alpha_0-2p_1q_2p_2,
\]
we f\/ind that
$c_{\infty,0}=c_{\infty,-1}=\cdots=c_{\infty,-(k-1)}=0$, $a_{\infty,n_0}b_{\infty,1}+c_{\infty,-k}d_{\infty,n_3}=0$.
Furthermore,
considering that
\[
tq_1^{\prime} =2q_1^2p_1-q_1^2+(\alpha_0+\alpha_1+\alpha_3)q_1-t+4tp_2+2q_1q_2p_2,
\]
we can show the lemma.
\end{proof}

\begin{proposition}
For $A_5^{(2)}(\alpha_j)_{0\leq j \leq 3}$,
there exists no solution such that
$q_1$, $p_1$, $p_2$ all have a~pole at $t=\infty$ and $q_2$ is holomorphic at $t=\infty$.
\end{proposition}

\begin{proof}
We treat the case where $n_1=1$. The other cases can be proved in the same way.
Comparing the coef\/f\/icients of the terms $t^{n_0+1}$ in
\[
tp_1^{\prime} =-2q_1p_1^2+2q_1p_1-(\alpha_0+\alpha_1+\alpha_3)p_1+\alpha_0-2p_1q_2p_2,
\]
we have
\[
-2a_{\infty,n_0}b_{\infty,0}-2a_{\infty,n_0-1}b_{\infty,1}+2a_{\infty,n_0}-2c_{\infty,-k}d_{\infty,n_3-1}-2c_{\infty,-k-1}d_{\infty,n_3}=0.
\]
\par
If $n_0-k\geq 3$,
comparing the coef\/f\/icients of the terms $t^{2n_0}$ in
\[
tq_1^{\prime} =2q_1^2p_1-q_1^2+(\alpha_0+\alpha_1+\alpha_3)q_1-t+4tp_2+2q_1q_2p_2,
\]
we obtain
\[
2b_{\infty,1}a_{\infty,n_0-1}+2a_{\infty,n_0}b_{\infty,0}-a_{\infty,n_0}+2c_{\infty,-k}d_{\infty,n_3-1}+2c_{\infty,-k-1}d_{\infty,n_3}=0.
\]
Then, it follows that $a_{\infty,n_0}=0$, which is impossible.

If $n_0-k=2$,
comparing the coef\/f\/icients of the terms $t^2$, $t^{3n_0-1}$ in
\begin{gather*}
tq_2^{\prime} =2q_2^2p_2-q_2^2+(\alpha_0+\alpha_1+\alpha_3)q_2-t+4tp_1+2q_1p_1q_2, \\
tp_2^{\prime} =-2q_2p_2^2+2q_2p_2-(\alpha_0+\alpha_1+\alpha_3)p_2+\alpha_1-2q_1p_1p_2,
\end{gather*}
we have
\begin{gather*}
2d_{\infty,n_3}c_{\infty,-k}c_{\infty,-k-1}+2d_{\infty,n_3-1}c_{\infty,-k}^2+4b_{\infty,1}\\
\qquad{} +
2c_{\infty,-k}a_{\infty,n_0}b_{\infty,0}+
2c_{\infty,-k}a_{\infty,n_0-1}b_{\infty,1}=0,     \\
-2c_{\infty,-k}d_{\infty,n_3-1}-2c_{\infty,-k-1}d_{\infty,n_3}
-2a_{\infty,n_0}b_{\infty,0}-2a_{\infty,n_0-1}b_{\infty,1}=0,
\end{gather*}
respectively.
Then,
it follows that $b_{\infty,1}=0$, which is impossible.

If $n_0-k=1$,
comparing the coef\/f\/icients of the terms $t^2$ in
\begin{equation*}
tq_2^{\prime} =2q_2^2p_2-q_2^2+(\alpha_0+\alpha_1+\alpha_3)q_2-t+4tp_1+2q_1p_1q_2,
\end{equation*}
we obtain
\[
2c_{\infty,-k}^2d_{\infty,n_3}+4b_{\infty,1}+2a_{\infty,n_0}b_{\infty,1}c_{\infty,-k}=4b_{\infty,1}=0,
\]
which is impossible.
\end{proof}

\subsection[The case where all of $(q_1,p_1,q_2,p_2)$ have a pole at $t=\infty$]{The case where all of $\boldsymbol{(q_1,p_1,q_2,p_2)}$ have a pole at $\boldsymbol{t=\infty}$}

\begin{lemma}
Suppose that
for $A_5^{(2)}(\alpha_j)_{0\leq j \leq 3}$,
there exists a solution such that
$q_1$, $p_1$, $q_2$, $p_2$ all have a pole at $t=\infty$.
Moreover,
assume that
$q_1$, $p_1$, $q_2$, $p_2$ have a pole of order $n_0$, $n_1$, $n_2$, $n_3$ $(n_0,n_1,n_2,n_3\geq 1)$ at $t=\infty$,
respectively.
Then,
$n_0+n_1=n_2+n_3$.
\end{lemma}

\begin{proof}
Considering
\[
tp_1^{\prime} =-2q_1p_1^2+2q_1p_1-(\alpha_0+\alpha_1+\alpha_3)p_1+\alpha_0-2p_1q_2p_2,
\]
we can show the lemma.
\end{proof}

Therefore,
we see that $n_0+n_1=n_2+n_3\geq 2$.

\begin{proposition}
For $A_5^{(2)}(\alpha_j)_{0\leq j \leq 3}$,
there exists no solution such that
$q_1$, $p_1$, $q_2$, $p_2$ all have a pole at $t=\infty$.
\end{proposition}

\begin{proof}
We treat the case where $n_0+n_1=n_2+n_3=2$.
The other cases can be proved in the same way.

Comparing the coef\/f\/icients of the term $t^3$ in
\[
tp_1^{\prime} =-2q_1p_1^2+2q_1p_1-(\alpha_0+\alpha_1+\alpha_3)p_1+\alpha_0-2p_1q_2p_2,
\]
we have $a_{\infty,1}b_{\infty,1}+c_{\infty,1}d_{\infty,1}=0$.

Comparing the coef\/f\/icients of the term $t^2$ in
\begin{gather*}
tq_1^{\prime} =2q_1^2p_1-q_1^2+(\alpha_0+\alpha_1+\alpha_3)q_1-t+4tp_2+2q_1q_2p_2,  \\
tp_1^{\prime} =-2q_1p_1^2+2q_1p_1-(\alpha_0+\alpha_1+\alpha_3)p_1+\alpha_0-2p_1q_2p_2,  \\
tq_2^{\prime} =2q_2^2p_2-q_2^2+(\alpha_0+\alpha_1+\alpha_3)q_2-t+4tp_1+2q_1p_1q_2, \\
tp_2^{\prime} =-2q_2p_2^2+2q_2p_2-(\alpha_0+\alpha_1+\alpha_3)p_2+\alpha_1-2q_1p_1p_2,
\end{gather*}
we obtain{\samepage
\begin{gather}
2a_{\infty,1}a_{\infty,0}b_{\infty,1}+2b_{\infty,0}a_{\infty,1}^2-a_{\infty,1}^2+4d_{\infty,1}
+2a_{\infty,1}c_{\infty,1}d_{\infty,0}+2a_{\infty,1}c_{\infty,0}d_{\infty,1}=0,   \nonumber \\
-2a_{\infty,1}b_{\infty,0}-2a_{\infty,0}b_{\infty,1}+2a_{\infty,1}-2c_{\infty,1}d_{\infty,0}-2c_{\infty,0}d_{\infty,1}=0,   \nonumber\\
2c_{\infty,1}c_{\infty,0}d_{\infty,1}+2d_{\infty,0}c_{\infty,1}^2-c_{\infty,1}^2+4b_{\infty,1}
+2c_{\infty,1}a_{\infty,1}b_{\infty,0}+2c_{\infty,1}a_{\infty,0}b_{\infty,1}=0,   \nonumber \\
-2c_{\infty,1}d_{\infty,0}-2c_{\infty,0}d_{\infty,1}+2c_{\infty,1}-2a_{\infty,1}b_{\infty,0}-2a_{\infty,0}b_{\infty,1}=0,  \label{*1}
\end{gather}}

\noindent
respectively.
Based on the second and fourth equations of \eqref{*1},
we have $a_{\infty,1}=c_{\infty,1}$.
From the f\/irst and second equations of \eqref{*1},
we obtain $a_{\infty,1}^2+4d_{\infty,1}=0$.
From the third and fourth equations of \eqref{*1},
we have $c_{\infty,1}^2+4b_{\infty,1}=0$.

Therefore,
since $a_{\infty,1}b_{\infty,1}+c_{\infty,1}d_{\infty,1}=0$,
it follows that
$a_{\infty,1}b_{\infty,1}=0$,
which is impos\-sible.
\end{proof}

\subsection{Summary}

\begin{proposition}
\label{prop:mero-solu-t=inf}
For $A_5^{(2)}(\alpha_j)_{0\leq j \leq 3}$,
there exists a meromorphic solution at $t=\infty$.
Then,
$q_1$, $p_1$, $q_2$, $p_2$ are uniquely expanded as follows:
\begin{gather*}
\begin{cases}
q_1=(-2\alpha_0+\alpha_3)   +\cdots, \\
p_1=1/4   +(-2\alpha_1+\alpha_3)(-2\alpha_1-\alpha_3) t^{-1}/4+\cdots, \\
q_2=(-2\alpha_1+\alpha_3)   +\cdots, \\
p_2=1/4   +(-2\alpha_0+\alpha_3)(-2\alpha_0-\alpha_3)t^{-1}/4 +\cdots.
\end{cases}
\end{gather*}
\end{proposition}

\section[Meromorphic solution at $t=0$]{Meromorphic solution at $\boldsymbol{t=0}$}\label{section2}

In this section,
we treat meromorphic solutions at $t=0$.
Then,
in the same way as Proposition~\ref{prop:mero-solu-t=inf},
we can show the following proposition:

\begin{proposition}
Suppose that
for
$A_5^{(2)}(\alpha_j)_{0\leq j \leq 3}$,
there exists a meromorphic solution at $t=0$.
Then, one of the following occurs:
\begin{enumerate}\itemsep=0pt
\item[$(1)$] $q_1$, $p_1$, $q_2$, $p_2$ are all holomorphic at $t=0$,

\item[$(2)$] $p_1$ has a pole of order one at $t=0$ and $q_1$, $q_2$, $p_2$ are all holomorphic at $t=0$,

\item[$(3)$] $p_2$ has a pole of order one at $t=0$ and $q_1$, $p_1$, $q_2$ are all holomorphic at $t=0$.
\end{enumerate}
\end{proposition}

In this paper,
we def\/ine the coef\/f\/icients of the Lauren series of $q_1$, $p_1$, $q_2$, $p_2$ at $t=0$
by
$a_{0,k}$, $b_{0,k}$, $c_{0,k}$, $d_{0,k}$, $k\in\mathbb{Z}$.
In this section,
we prove that the constant terms of $q_1$, $q_2$ at $t=0$,
$a_{0,0}$, $c_{0,0}$ are zero, or expressed by the parameters, $\alpha_j$ $(0\leq j \leq 3)$.

\subsection[The case where $q_1$, $p_1$, $q_2$, $p_2$ are all holomorphic at $t=0$]{The case where $\boldsymbol{q_1}$, $\boldsymbol{p_1}$, $\boldsymbol{q_2}$, $\boldsymbol{p_2}$ are all holomorphic at $\boldsymbol{t=0}$}

\begin{proposition}
\label{prop:t=0-holo}
Suppose that
for
$A_5^{(2)}(\alpha_j)_{0\leq j \leq 3}$,
there exists a solution such that
$q_1$, $p_1$, $q_2$, $p_2$ are all holomorphic at $t=0$.
Then, one of the following occurs:
\begin{enumerate}\itemsep=0pt
\item[$(1)$]
$a_{0,0}=0$, $-(\alpha_0+\alpha_1+\alpha_3)b_{0,0}+\alpha_0=0$, $c_{0,0}=0$, $-(\alpha_0+\alpha_1+\alpha_3)d_{0,0}+\alpha_1=0$,
\item[$(2)$]
$a_{0,0}=0$, $(-\alpha_0+\alpha_1-\alpha_3)b_{0,0}+\alpha_0=0$, $c_{0,0}=\alpha_0-\alpha_1+\alpha_3$, $(-\alpha_0+\alpha_1-\alpha_3)d_{0,0}-\alpha_1=0$,
\item[$(3)$]
$a_{0,0}=-\alpha_0+\alpha_1+\alpha_3$, $(\alpha_0-\alpha_1-\alpha_3)b_{0,0}-\alpha_0=0$, $c_{0,0}=0$, $(\alpha_0-\alpha_1-\alpha_3)d_{0,0}+\alpha_1=0$,
\item[$(4)$]
$a_{0,0}=-\alpha_0-\alpha_1+\alpha_3$, $(\alpha_0+\alpha_1-\alpha_3)b_{0,0}-\alpha_0=0$, $c_{0,0}=-\alpha_0-\alpha_1+\alpha_3$,
 $(\alpha_0+\alpha_1-\alpha_3)d_{0,0}-\alpha_1=0$.
 \end{enumerate}
\end{proposition}

\subsection[The case where $p_1$ has a pole at $t=0$]{The case where $\boldsymbol{p_1}$ has a pole at $\boldsymbol{t=0}$}

\begin{proposition}
\label{prop:t=0-p_1}
Suppose that
for
$A_5^{(2)}(\alpha_j)_{0\leq j \leq 3}$,
there exists a solution such that
$p_1$ has a pole at $t=0$
and $q_1$, $q_2$, $p_2$ are all holomorphic at $t=0$.
Then,
\begin{gather*}
\begin{cases}
q_1=(-8\alpha_0-8\alpha_3+6)t/\{(4\alpha_1-1)(4\alpha_1+1)\}+\cdots,   \\
p_1=(4\alpha_1-1)(4\alpha_1+1)t^{-1}/16+\cdots, \\
q_2=(-2\alpha_1+1/2)+\cdots, \\
p_2=1/4+\cdots.
\end{cases}
\end{gather*}
\end{proposition}

\subsection[The case where $p_2$ has a pole at $t=0$]{The case where $\boldsymbol{p_2}$ has a pole at $\boldsymbol{t=0}$}

\begin{proposition}
\label{prop:t=0-p_2}
Suppose that
for
$A_5^{(2)}(\alpha_j)_{0\leq j \leq 3}$,
there exists a solution such that
$p_2$ has a pole at $t=0$
and $q_1$, $p_1$, $q_2$ are all holomorphic at $t=0$.
Then,
\begin{gather*}
\begin{cases}
q_1=(-2\alpha_0+1/2)+\cdots, \\
p_1=1/4+\cdots, \\
q_2=(-8\alpha_1-8\alpha_3+6)t/\{(4\alpha_0-1)(4\alpha_0+1)\}+\cdots,   \\
p_2=(4\alpha_0-1)(4\alpha_0+1)t^{-1}/16+\cdots.
\end{cases}
\end{gather*}
\end{proposition}

\section[Meromorphic solution at $t=c\in\mathbb{C}^{*}$]{Meromorphic solution at $\boldsymbol{t=c\in\mathbb{C}^{*}}$}\label{section3}

In this section,
we deal with meromorphic solutions at $t=c\in\mathbb{C}^{*}$,
where $\mathbb{C}^{*}$ means the set of nonzero complex numbers.

\begin{proposition}
Suppose that
for
$A_5^{(2)}(\alpha_j)_{0\leq j \leq 3}$,
there exists a meromorphic solution at $t=c\in\mathbb{C}^{*}$
such that some of $(q_1,p_1,q_2,p_2)$ have a pole at $t=c$.
Then, one of the following occurs:
\begin{enumerate}\itemsep=0pt
\item[$(1)$] $q_1$ has a pole at $t=c$ and
$p_1,q_2,p_2$ are all holomorphic at $t=c$,
\item[$(2)$]  $q_2$ has a pole at $t=c$ and
$q_1,p_1,p_2$ are all holomorphic at $t=c$,
\item[$(3)$]  $q_1$, $q_2$ have both a pole at $t=c$ and
$p_1,p_2$ are both holomorphic at $t=c$,
\item[$(4)$] $q_1$, $p_2$ have both a pole at $t=c$ and
$p_1$, $q_2$ are both holomorphic at $t=c$,
\item[$(5)$]  $p_1$, $q_2$ have both a pole at $t=c$ and
$q_1$, $p_2$ are both holomorphic at $t=c$,
\item[$(6)$]  $p_1$, $p_2$ have both a pole at $t=c$ and
$q_1$, $q_2$ are both holomorphic at $t=c$,
\item[$(7)$]  $q_1$, $p_1$, $q_2$, $p_2$ all have a pole at $t=c$.
\end{enumerate}
\end{proposition}

\subsection[The case where $q_1$ has a pole at $t=c\in\mathbb{C}^{*}$]{The case where $\boldsymbol{q_1}$ has a pole at $\boldsymbol{t=c\in\mathbb{C}^{*}}$}

\begin{proposition}
Suppose that
for $A_5^{(2)}(\alpha_j)_{0\leq j \leq 3}$,
there exists a solution such that
$q_1$ has a pole at $t=c\in\mathbb{C}^{*}$ and
$p_1$, $q_2$, $p_2$ are all holomorphic at $t=c$.
Then,
either of the following occurs:
\begin{equation*}
(1) \
\begin{cases}
q_1=c(t-c)^{-1}+\cdots, \\
p_1=\displaystyle -\frac{\alpha_0}{c}(t-c)+\cdots, \\
\end{cases}\quad
(2) \
\begin{cases}
q_1=-c(t-c)^{-1}+\cdots, \\
p_1=\displaystyle 1+\frac{\alpha_1+\alpha_3}{c}(t-c)+\cdots,  \\
q_2=O(t-c),  \\
p_2=O(t-c).
\end{cases}
\end{equation*}
\end{proposition}

\subsection[The case where $q_2$ has a pole at $t=c\in\mathbb{C}^{*}$]{The case where $\boldsymbol{q_2}$ has a pole at $\boldsymbol{t=c\in\mathbb{C}^{*}}$}

\begin{proposition}
Suppose that
for $A_5^{(2)}(\alpha_j)_{0\leq j \leq 3}$,
there exists a solution such that
$q_2$ has a pole at $t=c\in\mathbb{C}^{*}$ and
$q_1$, $p_1$, $p_2$ are all holomorphic at $t=c$.
Then,
either of the following occurs:
\begin{equation*}
(1) \
\begin{cases}
q_2=c(t-c)^{-1}+\cdots, \\
p_2=\displaystyle -\frac{\alpha_1}{c}(t-c)+\cdots, \\
\end{cases}\quad
(2) \
\begin{cases}
q_1=O(t-c), \\
p_1=O(t-c),\\
q_2=-c(t-c)^{-1}+\cdots,  \\
p_2=\displaystyle 1+\frac{\alpha_0+\alpha_3}{c}(t-c)+\cdots.
\end{cases}
\end{equation*}
\end{proposition}

\subsection[The case where $q_1$, $q_2$ have a pole at $t=c\in\mathbb{C}^{*}$]{The case where $\boldsymbol{q_1}$, $\boldsymbol{q_2}$ have a pole at $\boldsymbol{t=c\in\mathbb{C}^{*}}$}

\begin{proposition}
Suppose that
for $A_5^{(2)}(\alpha_j)_{0\leq j \leq 3}$,
there exists a solution such that
$q_1$, $q_2$ have both a pole at $t=c\in\mathbb{C}^{*}$ and
$p_1$, $p_2$ are both holomorphic at $t=c$.
Then,
either of the following occurs:
\begin{equation*}
(1) \
\begin{cases}
q_1=-c(t-c)^{-1}+\cdots, \\
p_1=b_{c,0}+b_{c,1}(t-c)+\cdots, \\
q_2=-c(t-c)^{-1}+\cdots, \\
p_2=d_{c,0}+d_{c,1}(t-c)+\cdots,
\end{cases} \quad
(2) \
\begin{cases}
q_1=c(t-c)^{-1}+\cdots, \\
p_1=\displaystyle -\frac{\alpha_0}{c}(t-c)+\cdots, \\
q_2=c(t-c)^{-1}+\cdots, \\
p_2=\displaystyle -\frac{\alpha_1}{c}(t-c)+\cdots,
\end{cases}
\end{equation*}
where
$b_{c,0}+d_{c,0}=1$ and $b_{c,1}+d_{c,1}=  \frac{\alpha_3}{c}$.
\end{proposition}

\subsection[The case where $q_1$, $p_2$ have a pole at $t=c\in\mathbb{C}^{*}$]{The case where $\boldsymbol{q_1}$, $\boldsymbol{p_2}$ have a pole at $\boldsymbol{t=c\in\mathbb{C}^{*}}$}

\begin{proposition}
\label{prop:t=c(q_1,p_2)}
Suppose that
for $A_5^{(2)}(\alpha_j)_{0\leq j \leq 3}$,
there exists a solution such that
$q_1$, $p_2$ have both a pole at $t=c\in\mathbb{C}^{*}$ and
$p_1$, $q_2$ are both holomorphic at $t=c$.
Then, one of the following occurs:
\begin{gather*}
(1) \
\begin{cases}
q_1=4c(t-c)^{-1}+8/3+\cdots,     \\
p_1=0-\alpha_0/\{5c\}\cdot(t-c)+\cdots,       \\
q_2=(t-c)+(3\alpha_0-\alpha_1-\alpha_3+2)/\{2c\} \cdot (t-c)^2+\cdots,    \\
p_2=c(t-c)^{-2}-4\alpha_0/5 \cdot (t-c)^{-1}+\cdots,
\end{cases}
\\
(2) \
\begin{cases}
q_1=-c(t-c)^{-1}+(-1/4-\alpha_0)+\cdots,    \\
p_1=0+\alpha_0/\{5c\}\cdot(t-c)       \\
q_2=(t-c)+(-3\alpha_0-\alpha_1-\alpha_3+2)/\{2c\} \cdot (t-c)^2+\cdots,     \\
p_2=c(t-c)^{-2}+4\alpha_0/5 \cdot (t-c)^{-1}+\cdots,
\end{cases}
\\
(3) \
\begin{cases}
q_1=c(t-c)^{-1}+(3/4-\alpha_0 )+\cdots,     \\
p_1=1/2-1/\{12 c\} \cdot (t-c)+\cdots,        \\
q_2=-(t-c)+[(\alpha_1+\alpha_3)/c-3/\{4c\} ](t-c)^2+\cdots,      \\
p_2=-c/2\cdot (t-c)^{-2}-1/6(t-c)^{-1}+\cdots.
\end{cases}
\end{gather*}
\end{proposition}

\subsection[The case where $p_1$, $q_2$ have a pole at $t=c\in\mathbb{C}^{*}$]{The case where $\boldsymbol{p_1}$, $\boldsymbol{q_2}$ have a pole at $\boldsymbol{t=c\in\mathbb{C}^{*}}$}

\begin{proposition}
\label{prop:t=c(q_2,p_1)}
Suppose that
for $A_5^{(2)}(\alpha_j)_{0\leq j \leq 3}$,
there exists a solution such that
$p_1$, $q_2$ have both a pole at $t=c\in\mathbb{C}^{*}$ and
$q_1$, $p_2$ are both holomorphic at $t=c$.
Then, one of the following occurs:
\begin{gather*}
(1) \
\begin{cases}
q_1=(t-c)+(3\alpha_1-\alpha_0-\alpha_3+2)/\{2c\}\cdot (t-c)^2+\cdots, \\
p_1=c(t-c)^{-2}-4\alpha_1/5\cdot (t-c)^{-1}+\cdots, \\
q_2=4c(t-c)^{-1}+8/3+\cdots, \\
p_2=0-\alpha_1/\{5c\}\cdot (t-c)+\cdots,
\end{cases}
\\
(2) \
\begin{cases}
q_1=(t-c)+(-3\alpha_1-\alpha_0-\alpha_3+2)/\{2c\}\cdot (t-c)^2+\cdots,  \\
p_1=c(t-c)^{-2}+4\alpha_1/5\cdot (t-c)^{-1}+\cdots,  \\
q_2=-c(t-c)^{-1}+(-1/4-\alpha_1)+\cdots,  \\
p_2=0+\alpha_1/\{5c\}\cdot (t-c)+\cdots,
\end{cases}
\\
(3) \
\begin{cases}
q_1=-(t-c)+[(\alpha_0+\alpha_3)/c-3/\{4c\}](t-c)^2+\cdots,   \\
p_1=-c/2\cdot (t-c)^{-2}-1/6(t-c)^{-1}+\cdots,  \\
q_2=c(t-c)^{-1}+(3/4-\alpha_1)+\cdots,   \\
p_2=1/2-1/\{12c\}\cdot (t-c)+\cdots.
\end{cases}
\end{gather*}
\end{proposition}

\subsection[The case where $p_1$, $p_2$ have a pole at $t=c\in\mathbb{C}^{*}$]{The case where $\boldsymbol{p_1}$, $\boldsymbol{p_2}$ have a pole at $\boldsymbol{t=c\in\mathbb{C}^{*}}$}

\begin{proposition}
\label{prop:t=c(p_1,p_2)}
Suppose that
for $A_5^{(2)}(\alpha_j)_{0\leq j \leq 3}$,
there exists a solution such that
$p_1$, $p_2$ have both a pole at $t=c\in\mathbb{C}^{*}$ and
$q_1$, $q_2$ are both holomorphic at $t=c$.
Then,
\begin{gather*}
\begin{cases}
q_1=(-4d_{c,-1})+a_{c,1}(t-c)+\cdots,   \\
p_1=b_{c,-1}(t-c)^{-1}+(3/8+2b_{c,-1}^2/c)+\cdots, \\
q_2=(-4b_{c,-1})+c_{c,1}(t-c)+\cdots,   \\
p_2=d_{c,-1}(t-c)^{-1}+(3/8+2d_{c,-1}^2/c)+\cdots,
\end{cases}
\end{gather*}
where the coefficients satisfy
\begin{gather*}
16b_{c,-1}d_{c,-1}+c=0, \qquad
a_{c,1}b_{c,-1}c+c_{c,1}d_{c,-1}c+\frac{c}{2}(\alpha_0+\alpha_1+\alpha_3)=\frac{c}{2}.
\end{gather*}
\end{proposition}

\subsection[The case where $q_1$, $p_1$, $q_2$, $p_2$ have a pole at $t=c\in\mathbb{C}^{*}$]{The case where $\boldsymbol{q_1}$, $\boldsymbol{p_1}$, $\boldsymbol{q_2}$, $\boldsymbol{p_2}$ have a pole at $\boldsymbol{t=c\in\mathbb{C}^{*}}$}

\begin{proposition}
\label{prop:t=c-all}
Suppose that
for $A_5^{(2)}(\alpha_j)_{0\leq j \leq 3}$,
there exists a solution such that
$q_1$, $p_1$, $q_2$, $p_2$ all have a pole at $t=c\in\mathbb{C}^{*}$.
Then,
\begin{gather*}
\begin{cases}
q_1=-2c(t-c)^{-1}+(\sqrt{c}-4/3)+a_{c,1}(t-c)+\cdots,  \\
p_1=\sqrt{c}/4\cdot(t-c)^{-1}+1/2+b_{c,1}(t-c)+\cdots, \\
q_2=-2c(t-c)^{-1}+(-\sqrt{c}-4/3)+c_{c,1}(t-c)+\cdots,  \\
p_2=-\sqrt{c}/4\cdot (t-c)^{-1}+1/2+d_{c,1}(t-c)+\cdots,
\end{cases}
\end{gather*}
where the coefficients satisfy
\begin{gather*}
b_{c,1}+d_{c,1}=\alpha_3/\{2c\}, \qquad
a_{c,1}\sqrt{c}-c_{c,1}\sqrt{c}=2+2\alpha_3-2\alpha_0-2\alpha_1.
\end{gather*}
\end{proposition}

\subsection{Summary}

\begin{proposition}\label{prop:t=c-summary} \quad
\begin{enumerate}\itemsep=0pt
\item[$(1)$]
Suppose that
for $A_5^{(2)}(\alpha_j)_{0\leq j \leq 3}$,
there exists a meromorphic solution at $t=c\in \mathbb{C}^{*}$.
Then,
$q_1$, $q_2$ have both a pole of order at most one at $t=c$
and
the residues of $q_1$, $q_2$ at $t=c$ are expressed by $nc$ $(n\in\mathbb{Z})$.

\item[$(2)$] Suppose that
for $A_5^{(2)}(\alpha_j)_{0\leq j \leq 3}$,
there exists a rational solution.
Then,
$a_{\infty,0}-a_{0,0}\in\mathbb{Z}$, $c_{\infty,0}-c_{0,0}\in\mathbb{Z}$.
\end{enumerate}
\end{proposition}

\begin{proof}
Case (1) is obvious.
Let us prove case (2).
From the discussions in Sections~\ref{section1},~\ref{section2} and~\ref{section3},
it follows that
\begin{equation*}
q_1 =  a_{\infty,0}+\sum_{j=1}^{m_1} \frac{n_{j} c_{j}}{t-c_{j}}, \qquad
q_2 =  c_{\infty,0}+\sum_{k=1}^{m_2} \frac{n^{\prime}_{k} c^{\prime}_{k}}{t-c^{\prime}_{k}},   \qquad n_j, n^{\prime}_k \in\mathbb{Z},
\end{equation*}
where $m_1$, $m_2$ are both positive integers
and
$c_k \in\mathbb{C}^{*}$ $(1\leq k \leq m_1)$ and $c^{\prime}_j \in\mathbb{C}^{*}$ $(1\leq j \leq m_2)$ are
poles of $q_1$ and $q_2$, respectively.
If $q_1$ or $q_2$ is holomorphic in $\mathbb{C}^{*}$,
then its second sum is considered to be zero.

Considering the constant terms of the Taylor series of $q_1$, $q_2$ at $t=0$,
we can prove the proposition.
\end{proof}

\section[The Laurent series of the Hamiltonian $H$]{The Laurent series of the Hamiltonian $\boldsymbol{H}$}\label{section4}

In this section,
for a meromorphic solution at $t=\infty, 0$,
we f\/irst compute the constant terms $h_{\infty,0}$, $h_{0,0}$ of the Laurent series of the Hamiltonian $H$ at $t=\infty,0$.
Moreover,
for a meromorphic solution at $t=c\in\mathbb{C}^{*}$,
we calculate the residue of $H$ at $t=c$.

\subsection[The Laurent series of $H$ at $t=\infty$]{The Laurent series of $\boldsymbol{H}$ at $\boldsymbol{t=\infty}$}

\begin{proposition}
Suppose that
for $A_5^{(2)}(\alpha_j)_{0\leq j \leq 3}$,
there exists a meromorphic solution at $t=\infty$.
Then,
\[
h_{\infty,0}=\frac34(\alpha_0+\alpha_1+\alpha_3)^2
-\frac12(-2\alpha_0+\alpha_3)(-2\alpha_1+\alpha_3)
-3(\alpha_0+\alpha_1)\alpha_3.
\]
\end{proposition}

\subsection[The Laurent series of $H$ at $t=0$]{The Laurent series of $\boldsymbol{H}$ at $\boldsymbol{t=0}$}

\subsubsection[The case where $q_1$, $p_1$, $q_2$, $p_2$ are all holomorphic at $t=0$]{The case where $\boldsymbol{q_1}$, $\boldsymbol{p_1}$, $\boldsymbol{q_2}$, $\boldsymbol{p_2}$ are all holomorphic at $\boldsymbol{t=0}$}

\begin{proposition}
Suppose that
for $A_5^{(2)}(\alpha_j)_{0\leq j \leq 3}$,
there exists a solution
such that $q_1$, $p_1$, $q_2$, $p_2$ are all holomorphic at $t=0$.
Then,
\begin{equation*}
h_{0,0}=
\begin{cases}
0 & \text{if   case  $(1)$   occurs   in   Proposition~{\rm \ref{prop:t=0-holo}}}, \\
-\alpha_1(\alpha_0+\alpha_3) & \text{if   case   $(2)$   occurs   in   Proposition~{\rm \ref{prop:t=0-holo}}},  \\
-\alpha_0(\alpha_1+\alpha_3) & \text{if   case  $(3)$  occurs   in   Proposition~{\rm \ref{prop:t=0-holo}}},  \\
-\alpha_3(\alpha_0+\alpha_1) & \text{if   case   $(4)$   occurs  in  Proposition~{\rm \ref{prop:t=0-holo}}}.
\end{cases}
\end{equation*}
\end{proposition}

\subsubsection[The case where $p_1$ has a pole at $t=0$]{The case where $\boldsymbol{p_1}$ has a pole at $\boldsymbol{t=0}$}

\begin{proposition}
Suppose that
for $A_5^{(2)}(\alpha_j)_{0\leq j \leq 3}$,
there exists a solution
such that
$p_1$ has a pole at $t=0$
and
$q_1$, $q_2$, $p_2$ are all holomorphic at $t=0$.
Then,
\[
h_{0,0}=-\frac14(\alpha_0+\alpha_1+\alpha_3)^2+\alpha_1^2+\frac{3}{16}.
\]
\end{proposition}

\subsubsection[The case where $p_2$ has a pole at $t=0$]{The case where $\boldsymbol{p_2}$ has a pole at $\boldsymbol{t=0}$}

\begin{proposition}
Suppose that
for $A_5^{(2)}(\alpha_j)_{0\leq j \leq 3}$,
there exists a solution
such that
$p_2$ has a pole at $t=0$
and
$q_1$, $p_1$, $q_2$ are all holomorphic at $t=0$.
Then,
\[
h_{0,0}=-\frac14(\alpha_0+\alpha_1+\alpha_3)^2+\alpha_0^2+\frac{3}{16}.
\]
\end{proposition}

\subsection[The Laurent series of $H$ at $t=c\in\mathbb{C}^{*}$]{The Laurent series of $\boldsymbol{H}$ at $\boldsymbol{t=c\in\mathbb{C}^{*}}$}

\subsubsection[The case where $q_1$ has a pole at $t=c\in\mathbb{C}^{*}$]{The case where $\boldsymbol{q_1}$ has a pole at $\boldsymbol{t=c\in\mathbb{C}^{*}}$}

\begin{proposition}
Suppose that
for $A_5^{(2)}(\alpha_j)_{0\leq j \leq 3}$,
there exists a solution such that
$q_1$ has a pole at $t=c\in\mathbb{C}^{*}$ and
$p_1$, $q_2$, $p_2$ are all holomorphic at $t=c$.
Then,
$H$ is holomorphic at $t=c$.
\end{proposition}

\subsubsection[The case where $q_2$ has a pole at $t=c\in\mathbb{C}^{*}$]{The case where $\boldsymbol{q_2}$ has a pole at $\boldsymbol{t=c\in\mathbb{C}^{*}}$}

\begin{proposition}
Suppose that
for $A_5^{(2)}(\alpha_j)_{0\leq j \leq 3}$,
there exists a solution such that
$q_2$ has a pole at $t=c\in\mathbb{C}^{*}$ and
$q_1$, $p_1$, $p_2$ are all holomorphic at $t=c$.
Then,
$H$ is holomorphic at $t=c$.
\end{proposition}

\subsubsection[The case where $q_1$, $q_2$ have a pole at $t=c\in\mathbb{C}^{*}$]{The case where $\boldsymbol{q_1}$, $\boldsymbol{q_2}$ have a pole at $\boldsymbol{t=c\in\mathbb{C}^{*}}$}

\begin{proposition}
Suppose that
for $A_5^{(2)}(\alpha_j)_{0\leq j \leq 3}$,
there exists a solution such that
$q_1$, $q_2$ have both a pole at $t=c\in\mathbb{C}^{*}$ and
$p_1$, $p_2$ are both holomorphic at $t=c$.
Then,
$H$ is holomorphic at $t=c$.
\end{proposition}

\subsubsection[The case where $q_1$, $p_2$ have a pole at $t=c\in\mathbb{C}^{*}$]{The case where $\boldsymbol{q_1}$, $\boldsymbol{p_2}$ have a pole at $\boldsymbol{t=c\in\mathbb{C}^{*}}$}

\begin{proposition}
Suppose that
for $A_5^{(2)}(\alpha_j)_{0\leq j \leq 3}$,
there exists a solution such that
$q_1$, $q_2$ have both a pole at $t=c\in\mathbb{C}^{*}$ and
$p_1$, $p_2$ are both holomorphic at $t=c$.
Then,
$H$ has a pole of order one at $t=c$
and
\begin{equation*}
\mathop{\mathrm{Res}}\limits_{t=c}H=
\begin{cases}
c & \text{if   case $(1)$ occurs   in   Proposition~{\rm \ref{prop:t=c(q_1,p_2)}}},   \\
c & \text{if   case   $(2)$   occurs   in  Proposition~{\rm \ref{prop:t=c(q_1,p_2)}}},   \\
c/2 & \text{if   case $(3)$ occurs   in   Proposition~{\rm \ref{prop:t=c(q_1,p_2)}}}.
\end{cases}
\end{equation*}
\end{proposition}

\subsubsection[The case where $p_1$, $q_2$ have a pole at $t=c\in\mathbb{C}^{*}$]{The case where $\boldsymbol{p_1}$, $\boldsymbol{q_2}$ have a pole at $\boldsymbol{t=c\in\mathbb{C}^{*}}$}

\begin{proposition}
Suppose that
for $A_5^{(2)}(\alpha_j)_{0\leq j \leq 3}$,
there exists a solution such that
$p_1$, $q_2$ have both a pole at $t=c\in\mathbb{C}^{*}$ and
$q_1$, $p_2$ are both holomorphic at $t=c$.
Then,
$H$ has a pole of order one at $t=c$
and
\begin{equation*}
\mathop{\mathrm{Res}}\limits_{t=c}H=
\begin{cases}
c & \text{if   case  $(1)$  occurs   in  Proposition~{\rm \ref{prop:t=c(q_2,p_1)}}},   \\
c & \text{if  case  $(2)$   occurs   in   Proposition~{\rm \ref{prop:t=c(q_2,p_1)}}},  \\
c/2 & \text{if    case   $(3)$   occurs   in   Proposition~{\rm \ref{prop:t=c(q_2,p_1)}}}.
\end{cases}
\end{equation*}
\end{proposition}

\subsubsection[The case where $p_1$, $p_2$ have a pole at $t=c\in\mathbb{C}^{*}$]{The case where $\boldsymbol{p_1}$, $\boldsymbol{p_2}$ have a pole at $\boldsymbol{t=c\in\mathbb{C}^{*}}$}

\begin{proposition}
Suppose that
for $A_5^{(2)}(\alpha_j)_{0\leq j \leq 3}$,
there exists a solution such that
$p_1$, $p_2$ have both a pole at $t=c\in\mathbb{C}^{*}$ and
$q_1$, $q_2$ are both holomorphic at $t=c$.
Then,
$H$ has a pole of order one at $t=c$
and
$\mathop{\mathrm{Res}}\limits_{t=c}H=c/4$.
\end{proposition}

\subsubsection[The case where $q_1$, $p_1$, $q_2$, $p_2$ have a pole at $t=c\in\mathbb{C}^{*}$]{The case where $\boldsymbol{q_1}$, $\boldsymbol{p_1}$, $\boldsymbol{q_2}$, $\boldsymbol{p_2}$ have a pole at $\boldsymbol{t=c\in\mathbb{C}^{*}}$}

\begin{proposition}
Suppose that
for $A_5^{(2)}(\alpha_j)_{0\leq j \leq 3}$,
there exists a solution such that
$q_1$, $p_1$, $q_2$, $p_2$ all have a pole at $t=c\in\mathbb{C}^{*}$.
Then,
$H$ has a pole of order one at $t=c$ and
$\mathop{\mathrm{Res}}\limits_{t=c} H=c/4$.
\end{proposition}

\subsection{Summary}

\begin{proposition}\quad
\begin{enumerate}\itemsep=0pt
\item[$(1)$] Suppose that
for $A_5^{(2)}(\alpha_j)_{0\leq j \leq 3}$,
there exists a meromorphic solution at $t=c\in\mathbb{C}^{*}$.
Then,
the residue of $H$ at $t=c$ is expressed by $nc/4$ $(n\in\mathbb{Z})$.

\item[$(2)$]
Suppose that
for $A_5^{(2)}(\alpha_j)_{0\leq j \leq 3}$,
there exists a rational solution.
Then,
$4(h_{\infty,0}-h_{0,0})\in\mathbb{Z}$.
\end{enumerate}
\end{proposition}

\begin{proof}
Case (1) is obvious.
Case (2) can be proved in the same way
as Proposition~\ref{prop:t=c-summary}.
\end{proof}

\section{Necessary condition \dots\ (1)}\label{section5}

\subsection[The case where $q_1$, $p_1$, $q_2$, $p_2$ are all holomorphic at $t=0$]{The case where $\boldsymbol{q_1}$, $\boldsymbol{p_1}$, $\boldsymbol{q_2}$, $\boldsymbol{p_2}$ are all holomorphic at $\boldsymbol{t=0}$}

\subsubsection[The case where $a_{0,0}=0$, $c_{0,0}=0$]{The case where $\boldsymbol{a_{0,0}=0}$, $\boldsymbol{c_{0,0}=0}$}

\begin{proposition}
\label{prop:necssary-holo-(1)}
Suppose that
for $A_5^{(2)}(\alpha_j)_{0\leq j \leq 3}$,
there exists a rational solution such that
$q_1$, $p_1$, $q_2$, $p_2$ are all holomorphic at $t=0$.
Moreover,
assuming that
$a_{0,0}=0$, $c_{0,0}=0$,
then,
$-2\alpha_0+\alpha_3\in\mathbb{Z}$, $-2\alpha_1+\alpha_3\in\mathbb{Z}$.
\end{proposition}

\begin{proof}
The proposition follows from Propositions \ref{prop:mero-solu-t=inf}, \ref{prop:t=c-summary}.
\end{proof}

\subsubsection[The case where $a_{0,0}=0$, $c_{0,0}\neq0$]{The case where $\boldsymbol{a_{0,0}=0}$, $\boldsymbol{c_{0,0}\neq0}$}

\begin{proposition}
\label{prop:necssary-holo-(2)}
Suppose that
for $A_5^{(2)}(\alpha_j)_{0\leq j \leq 3}$,
there exists a rational solution such that
$q_1$, $p_1$, $q_2$, $p_2$ are all holomorphic at $t=0$.
Moreover,
assuming that
$a_{0,0}=0$, $c_{0,0}\neq0$,
then,
$-2\alpha_0+\alpha_3\in\mathbb{Z}$, $2\alpha_1+\alpha_3\in\mathbb{Z}$.
\end{proposition}

\begin{proof}
The proposition follows from Propositions \ref{prop:mero-solu-t=inf}, \ref {prop:t=0-holo} and \ref{prop:t=c-summary}.
\end{proof}

\subsubsection[The case where $a_{0,0}\neq0$, $c_{0,0}=0$]{The case where $\boldsymbol{a_{0,0}\neq0}$, $\boldsymbol{c_{0,0}=0}$}

\begin{proposition}
\label{prop:necssary-holo-(3)}
Suppose that
for $A_5^{(2)}(\alpha_j)_{0\leq j \leq 3}$,
there exists a rational solution such that
$q_1$, $p_1$, $q_2$, $p_2$ are all holomorphic at $t=0$.
Moreover,
assuming that
$a_{0,0}\neq0$, $c_{0,0}=0$,
then,
$2\alpha_0+\alpha_3\in\mathbb{Z}$, $-2\alpha_1+\alpha_3\in\mathbb{Z}$.
\end{proposition}

\begin{proof}
The proposition follows from Propositions \ref{prop:mero-solu-t=inf}, \ref {prop:t=0-holo} and \ref{prop:t=c-summary}.
\end{proof}

\subsubsection[The case where $a_{0,0}\neq0$, $c_{0,0}\neq0$]{The case where $\boldsymbol{a_{0,0}\neq0}$, $\boldsymbol{c_{0,0}\neq0}$}

\begin{proposition}
\label{prop:necssary-holo-(4)}
Suppose that
for $A_5^{(2)}(\alpha_j)_{0\leq j \leq 3}$,
there exists a rational solution such that
$q_1$, $p_1$, $q_2$, $p_2$ are all holomorphic at $t=0$.
Moreover,
assuming that
$a_{0,0}\neq0$, $c_{0,0}\neq0$,
then,
$2\alpha_0+\alpha_3\in\mathbb{Z}$, $2\alpha_1+\alpha_3\in\mathbb{Z}$.
\end{proposition}

\begin{proof}
From Propositions \ref{prop:mero-solu-t=inf}, \ref {prop:t=0-holo} and \ref{prop:t=c-summary},
it follows that $\alpha_0-\alpha_1\in\mathbb{Z}$.

If $\alpha_0\neq 0$,
by Proposition~\ref {prop:t=0-holo},
we f\/ind that
$s_0(q_1,p_1,q_2,p_2)$ is a rational solution of $A_5^{(2)}(-\alpha_0$, $\alpha_1,\alpha_2+\alpha_0,\alpha_3)$
such that all of $s_0(q_1,p_1,q_2,p_2)$ are holomorphic at $t=0$ and $a_{0,0}=0$, \mbox{$c_{0,0}\neq0$}.
Then,
from Proposition~\ref{prop:necssary-holo-(2)},
we obtain the necessary condition.
If $\alpha_1\neq 0$,
by $s_1$ and Proposition~\ref{prop:necssary-holo-(3)},
we obtain the necessary condition in the same way.

If $\alpha_0=\alpha_1=0$ and $\alpha_2\neq 0$,
by Proposition~\ref {prop:t=0-holo},
we see that
$s_2(q_1,p_1,q_2,p_2)$ is a rational solution of $A_5^{(2)}(\alpha_2,\alpha_2,-\alpha_2,\alpha_3+2\alpha_2)$
such that all of $s_0(q_1,p_1,q_2,p_2)$ are holomorphic at $t=0$ and $a_{0,0}\neq0$, $c_{0,0}\neq0$.
Based on the above discussion,
considering that $\alpha_0+\alpha_1+2\alpha_2+\alpha_3=1/2$,
we can obtain the necessary condition.

The remaining case is that $\alpha_0=\alpha_1=\alpha_2=0$, $\alpha_3=1/2$.
We prove that
for $A_5^{(2)}(0,0,0,1/2)$,
there exists no rational solution such that
$q_1$, $p_1$, $q_2$, $p_2$ are all holomorphic at $t=0$ and $a_{0,0}\neq0$, $c_{0,0}\neq0$.
If there exists such a rational solution,
by Proposition~\ref {prop:t=0-holo},
we f\/ind that $b_{0,0}=d_{0,0}=0$.
Then,
$s_3(q_1,p_1,q_2,p_2)$ is a rational solution of $A_5^{(2)}(0,0,1/2,-1/2)$
such that all of $s_3(q_1,p_1,q_2,p_2)$ are holomorphic at $t=0$ and $a_{0,0}=c_{0,0}=0$.
Therefore,
it follows from Proposition~\ref{prop:necssary-holo-(2)}
that $-2\cdot0+(-1/2)\in\mathbb{Z}$,
which is impossible.
\end{proof}

\subsection[The case where $p_1$ has a pole at $t=0$]{The case where $\boldsymbol{p_1}$ has a pole at $\boldsymbol{t=0}$}

\begin{proposition}
Suppose that
for $A_5^{(2)}(\alpha_j)_{0\leq j \leq 3}$,
there exists a rational solution such that
$p_1$ has a pole at $t=0$
and
$q_1$, $q_2$, $p_2$ are all holomorphic at $t=0$.
Then,
$-2\alpha_0+\alpha_3\in\mathbb{Z}$, $\alpha_3-1/2\in\mathbb{Z}$.
\end{proposition}

\begin{proof}
The proposition follows from Propositions~\ref{prop:mero-solu-t=inf},~\ref{prop:t=0-p_1} and~\ref{prop:t=c-summary}.
\end{proof}

By $s_1s_2$,
we can prove the following corollary.

\begin{corollary}
Suppose that
for $A_5^{(2)}(\alpha_j)_{0\leq j \leq 3}$,
there exists a rational solution such that
$p_1$ has a pole at $t=0$
and
$q_1$, $q_2$, $p_2$ are all holomorphic at $t=0$.
Then,
by some B\"acklund transformations,
the parameters
can be transformed so that $-2\alpha_0+\alpha_3\in \mathbb{Z}$, $-2\alpha_1+\alpha_3\in\mathbb{Z}$.
\end{corollary}

\subsection[The case where $p_2$ has a pole at $t=0$]{The case where $\boldsymbol{p_2}$ has a pole at $\boldsymbol{t=0}$}

\begin{proposition}
Suppose that
for $A_5^{(2)}(\alpha_j)_{0\leq j \leq 3}$,
there exists a rational solution such that
$p_2$ has a pole at $t=0$
 and
$q_1$, $p_1$, $q_2$ are all holomorphic at $t=0$.
Then,
$-2\alpha_1+\alpha_3\in\mathbb{Z}$, $\alpha_3-1/2\in\mathbb{Z}$.
\end{proposition}

\begin{proof}
The proposition follows from Propositions \ref{prop:mero-solu-t=inf}, \ref{prop:t=0-p_2} and \ref{prop:t=c-summary}.
\end{proof}

By $s_0s_2$,
we can prove the following corollary.

\begin{corollary}
Suppose that
for $A_5^{(2)}(\alpha_j)_{0\leq j \leq 3}$,
there exists a rational solution such that
$p_2$ has a pole at $t=0$
and
$q_1$, $p_1$, $q_2$ are all holomorphic at $t=0$.
Then,
by some B\"acklund transformations,
the parameters
can be transformed so that $-2\alpha_0+\alpha_3\in \mathbb{Z}$, $-2\alpha_1+\alpha_3\in\mathbb{Z}$.
\end{corollary}

\subsection{Summary}

\begin{proposition}
\label{prop:necessary-1}
Suppose that
for $A_5^{(2)}(\alpha_j)_{0\leq j \leq 3}$,
there exists a rational solution.
Then,
one of the following occurs:
\begin{alignat*}{3}
&(1) &\quad  -2\alpha_0+\alpha_3&\in\mathbb{Z}, &\quad  -2\alpha_1+\alpha_3&\in\mathbb{Z},  \\
&(2) &  -2\alpha_0+\alpha_3&\in\mathbb{Z}, &  2\alpha_1+\alpha_3&\in\mathbb{Z},   \\
&(3) &   2\alpha_0+\alpha_3&\in\mathbb{Z},  & -2\alpha_1+\alpha_3&\in\mathbb{Z},  \\
&(4) &   2\alpha_0+\alpha_3&\in\mathbb{Z},  &  2\alpha_1+\alpha_3&\in\mathbb{Z},  \\
&(5) &  -2\alpha_0+\alpha_3&\in\mathbb{Z},  & \alpha_3-1/2&\in\mathbb{Z},\\
&(6) &  -2\alpha_1+\alpha_3&\in\mathbb{Z},  & \alpha_3-1/2&\in\mathbb{Z}.
\end{alignat*}
\end{proposition}

\begin{corollary}
\label{coro:necessary-2}
Suppose that
for $A_5^{(2)}(\alpha_j)_{0\leq j \leq 3}$,
there exists a rational solution.
Then,
by some B\"acklund transformations,
the parameters can be transformed so that
$-2\alpha_0+\alpha_3\in\mathbb{Z}$, \mbox{$-2\alpha_1+\alpha_3\in\mathbb{Z}$}.
\end{corollary}

\section{Necessary condition \dots\ (2)}\label{section6}

\subsection{Shift operators}

In order to transform the parameters to the standard form,
let us construct shift operators.

\begin{proposition}
Let the shift operators $T_0$, $T_1$, $T_2$ be defined by
\begin{equation*}
T_0=\pi s_2s_3 s_2 s_1 s_0, \qquad
T_1=s_0T_0s_0, \qquad
T_2=s_2T_0s_2,
\end{equation*}
respectively.
Then,
\begin{gather*}
T_0(\alpha_0,\alpha_1,\alpha_2,\alpha_3)  =(\alpha_0+1/2,\alpha_1+1/2,\alpha_2-1/2,\alpha_3),     \\
T_1(\alpha_0,\alpha_1,\alpha_2,\alpha_3)  =(\alpha_0-1/2,\alpha_1+1/2,\alpha_2,\alpha_3),     \\
T_2(\alpha_0,\alpha_1,\alpha_2,\alpha_3)  =(\alpha_0,\alpha_1,\alpha_2+1/2,\alpha_3-1),
\end{gather*}
respectively.
\end{proposition}


\subsection{The properties of B\"acklund transformations}

\begin{proposition}\quad
\begin{enumerate}\itemsep=0pt
\item[$(1)$] If $p_1\equiv0$ for $A_5^{(2)}(\alpha_j)_{0\leq j \leq 3}$,
then $\alpha_0=0$.
\item[$(2)$] If $p_2\equiv 0$ for $A_5^{(2)}(\alpha_j)_{0\leq j \leq 3}$,
then $\alpha_1=0$.
\item[$(3)$] If $q_1q_2+t\equiv 0$ for $A_5^{(2)}(\alpha_j)_{0\leq j \leq 3}$,
then $\alpha_2=0$.
\item[$(4)$] If $p_1+p_2-1\equiv 0$ for $A_5^{(2)}(\alpha_j)_{0\leq j \leq 3}$,
then $\alpha_3=0$.
\end{enumerate}
\end{proposition}


By this proposition,
we can consider $s_0$ as the identical transformation,
if $p_0\equiv 0$.
In the same way,
we consider
each of
$s_1$, $s_2$, $s_3$ as the identical transformation,
if $p_2\equiv 0$, or
if $q_1q_2+t\equiv 0$, or
if $p_1+p_2-1\equiv 0$,
respectively.

\subsection{Reduction of the parameters to the standard form}

By Corollary \ref{coro:necessary-2},
using $T_0$,
we can transform the parameters to $(\alpha_0,\alpha_1,\alpha_2,\alpha_3)=(\alpha_3/2,\alpha_3/2$, $\alpha_2,\alpha_3)$.

\begin{proposition}
\label{prop:necessary-3}
Suppose that
for $A_5^{(2)}(\alpha_j)_{0\leq j \leq 3}$,
there exists a rational solution.
Then,
by some B\"acklund transformations,
the parameters can be transformed so that
$-2\alpha_0+\alpha_3=0$, $-2\alpha_1+\alpha_3=0$.
\end{proposition}

\section{Classif\/ication of rational solutions}\label{section7}

\subsection[Rational solution of $A_5^{(2)}(\alpha_3/2,\alpha_3/2,\alpha_2,\alpha_3)$]{Rational solution of $\boldsymbol{A_5^{(2)}(\alpha_3/2,\alpha_3/2,\alpha_2,\alpha_3)}$}

\begin{proposition}
\label{prop:standard}
For $A_5^{(2)}(\alpha_3/2,\alpha_3/2,\alpha_2,\alpha_3)$,
there exists a rational solution and
$(q_1,p_1,q_2$, $p_2)=(0,1/4,0,1/4)$.
Moreover,
it is unique.
\end{proposition}

\begin{proof}
The proposition follows from
the direct calculation and
Proposition~\ref{prop:uniqueness}.
\end{proof}

\subsection{Proof of main theorem}

Let us prove our main theorem.

\begin{proof}
Suppose that for $A_5^{(2)}(\alpha_j)_{0\leq j \leq 3}$,
there exists a rational solution.
Then,
from Proposition~\ref{prop:necessary-1},
we f\/ind that
the parameters satisfy one of the conditions in the theorem.
Moreover,
from Proposition~\ref{prop:necessary-3},
we see that
the parameters
can be transformed so that $-2\alpha_0+\alpha_3=-2\alpha_1+\alpha_3=0$.

From Proposition~\ref{prop:standard},
it follows that
for $A_5^{(2)}(\alpha_3/2,\alpha_3/2,\alpha_2,\alpha_3)$,
there exists a unique rational solution
such that
$(q_1,p_1,q_2,p_2)=(0,1/4,0,1/4)$,
which proves the main theorem.
\end{proof}

\appendix
\section{Examples of rational solutions}\label{appendixA}

In this appendix,
we give examples of rational solutions of $A_5^{(2)}(\alpha_j)_{0\leq j \leq 3}$.
For the purpose,
we use the shift operators, $T_0$, $T_1$, $T_2$, and the seed rational solution,
\[
(q_1,p_1,q_2,p_2)=(0,1/4,0,1/4)  \qquad \mathrm{for} \quad A_5^{(2)}(\alpha_3/2,\alpha_3/2,\alpha_2,\alpha_3).
\]
Then,
we obtain the following examples of rational solutions:
\\
for $A_5^{(2)}(\alpha_3/2+1/2,\alpha_3/2+1/2,\alpha_2-1/2,\alpha_3)$,
\begin{gather*}
(q_1,p_1,q_2,p_2)=
\left(
-1,\frac14+\frac{\alpha_3}{2(t+4\alpha_3^2)}-\frac{4\alpha_3+1}{4(t+1)},
-1,\frac14-\frac{\alpha_3}{2(t+4\alpha_3^2)}+\frac{4\alpha_3+1}{4(t+1)},
\right);
\end{gather*}
for $A_5^{(2)}(\alpha_3/2-1/2,\alpha_3/2+1/2,\alpha_2,\alpha_3)$,
\begin{gather*}
q_1 =2\alpha_3-\cfrac{1}{1+
                    \cfrac{\alpha_3(1-2\alpha_3)}{t}}
    +\cfrac{-2\alpha_3+2}{1+
        \cfrac{2\alpha_3+1}{-t+\alpha_3(2\alpha_3+1)-
        \cfrac{1}{1+
        \cfrac{\alpha_3(-2\alpha_3+1)}{t}}}},   \\
p_1 =\frac14+\cfrac{2\alpha_3+1}{-4t+4\alpha_3(2\alpha_3+1)-
                 \cfrac{1}{1-
                 \cfrac{\alpha_3(2\alpha_3-1)}{t}}},   \\
q_2 =-\cfrac{1}{1-
        \cfrac{\alpha_3(2\alpha_3-1)}{t}},     \\
p_2 =\frac14+\frac{\alpha_3(2\alpha_3-1)}{t}-
                   \cfrac{2\alpha_3+1}{
                   \cfrac{-4}{1+
                   \cfrac{\alpha_3(1-2\alpha_3)}{t}} +
                   \cfrac{4t}{2\alpha_3-
                   \cfrac{1}{1+
                   \cfrac{\alpha_3(1-2\alpha_3)}{t}}}};
\end{gather*}
for $A_5^{(2)}(\alpha_3/2,\alpha_3/2,\alpha_2+1/2,\alpha_3-1)$,
\begin{gather*}
(q_1,p_1,q_2,p_2)=
\bigg(
-1,\frac14-\frac{2\alpha_3-1}{4\{t+(2\alpha_3-1)^2 \}}+\frac{4\alpha_3-3}{2(t+1)},
-1,\\
\hphantom{(q_1,p_1,q_2,p_2)= \bigg(}{} \frac14+\frac{2\alpha_3-1}{4\{t+(2\alpha_3-1)^2 \}}-\frac{4\alpha_3-3}{2(t+1)}
\bigg).
\end{gather*}

\subsection*{Acknowledgments}

The author wishes to express his sincere thanks to Professor Yousuke Ohyama.
In addition,
he is also indebted the referees for their useful comments.

\newpage

\pdfbookmark[1]{References}{ref}
\LastPageEnding

\end{document}